\documentclass[12pt]{iopart}

\usepackage{tcolorbox}

\usepackage{iopams}

\newtheorem{thm}{Theorem}
\newtheorem{Prop}{Proposition}
\newtheorem{Def}{Definition}
\newtheorem{ass}{Assumption}
\newtheorem{Lem}{Lemma}
\newtheorem{remark}{Remark}
\newtheorem{col}{Corollary}

\begin{document}

\title[Incomplete crossing and semi-topological horseshoes]{Incomplete crossing and semi-topological horseshoes}

\author{Junfeng Cheng$^1$ and Xiao-Song Yang$^2$\footnote{Author to whom any correspondence should be addressed.}}

\address{$^1$School of Mathematics and Statistics,\\
	Huazhong University of Science and Technology,\\
	Wuhan 430074, R. P. China}

\address{$^2$School of Mathematics and Statistics\\
	Huazhong University of Science and Technology\\
	Wuhan, 430074, R. P. China\\
	Hubei Key Laboratory of Engineering Modeling and Scientific Computing,\\
	Huazhong University of Science and Technology,\\
	Wuhan 430074, R. P. China}

\ead{$^1$ydtq19@163.com and $^2$yangxs@hust.edu.cn}
\vspace{10pt}

\begin{abstract}
This paper enriches the topological horseshoe theory using finite subshift theory in symbolic dynamical systems, and develops an elementary framework addressing incomplete crossing and semi-horseshoes. Two illustrative examples are provided: one from the perturbed Duffing system and another from a polynomial system proposed by Chen, demonstrating the prevalence of semi-horseshoes in chaotic systems. Moreover, the semi-topological horseshoe theory enhances the detection of chaos and improves the accuracy of topological entropy estimation.
\end{abstract}

\ams{37B40, 37B10, 34D45}

\noindent{\it Keywords}: topological horseshoe, subshift, Duffing system,  polynomial system

%
%
%
%
%

\section{Introduction}

The exploration of nonlinear dynamics throughout the 20th century has revealed chaos as one of its most interesting and complex phenomena. The seminal work by E. N. Lorenz in 1963 \cite{Lorenz1963} introduced the Lorenz system,  which is a fundamental model in the study of chaotic motion. However, it was not until W. Tucker's work in 2002 that the chaotic nature of the Lorenz system was rigorously proven \cite{Tucker2002}, which shows that grasping the complexities of chaotic systems is an arduous task.

When exploring chaos in dynamical systems, the topological horseshoe theory is an important tool that cannot be overlooked. The exploration of topological horseshoe theory began with the proposal of the discovery of the horseshoe \cite{Smale1967}, which later became known as the Smale horseshoe. Subsequently, Kennedy and Yorke introduced the formal concept of the topological horseshoe \cite{Kennedy2001-1, Kennedy2001-2}, which was further enhanced by Yang and Tang \cite{Yang2004, Yang2009} to address piecewise continuous maps. With this theorem, we can more conveniently use numerical methods, such as the Runge-Kutta method, to study chaos.

Following the work in \cite{Yang2004, Yang2009}, in order to address richer and more complex chaotic behavior, we present theory of incomplete crossings and semi-topological horseshoes with the aid of subshifts theory in symbolic dynamical systems.
This theory not only assists in detecting chaos but also provides a relatively accurate lower bound estimation of the topological entropy of dynamical system.

As illustrative applications, we investigate the Duffing system perturbed by a periodic term and a polynomial system proposed by Chen \cite{Chen1999}.
Different from the well-established results obtained via Melnikov's method \cite{Melnikov1963}, we explore the horseshoe structure of the Duffing system under specific parameter settings when the external force $\gamma$ takes certain values. It is found that even for $\gamma$ not small enough, for example $\gamma=1.184$, the system is still chaotic.
 Using our theory, we derive that the Poincaré map of the perturbed Duffing system satisfies the inequality :
\[ h(P_{\gamma=1.184}) \geq \ln \frac{1+\sqrt{5}}{2}, \]
where $h(P_{\gamma=1.184})$ denotes the topological entropy of its Poincaré map $P_{\gamma=1.184}$.

Another example arises from the polynomial system proposed by Chen. In \cite{Cheng2024}, we observed a compelling crossing structure in the Poincaré map. Applying our theory to this polynomial system, we deduce that $h(P)$ satisfies the inequality:
\[ h(P) \geq \ln (2+\sqrt{2}), \]
where $h(P)$ represents the topological entropy of the Poincaré map $P$ of the polynomial system.
We observe that the same results can be derived via mathematical induction, albeit with greater complexity in notation and computation.

The paper is structured as follows: Section \ref{THT} provides an overview of topological horseshoe theory. Section \ref{pre} introduces the theory of subshifts in symbolic dynamical systems, establishing concepts for later sections. In Section \ref{ITHT}, we present the main results of our study. Section \ref{TwoBlocks} demonstrates the application of our theory in horseshoes involving two blocks. Section \ref{Duffing} and Section \ref{Chen} revisit the perturbed Duffing system and the Chen system using our theory. Finally, Section \ref{summary} concludes with a concise summary of our research findings.

\section{Topological Horseshoe}\label{THT}

For the purpose of this paper, we offer a review of the topological horseshoes in this section. 
Prior to describing the concept of a topological horseshoe, it is essential to introduce the definition of semi-conjugation. 
The following two definitions can be found in \cite{Wiggins1988}.

\begin{Def} 
	Let $M$ and $N$ be topological spaces and consider two continuous maps $f:M\to M$ and $g:N\to N$. The map $f$ is said to be semi-conjugate to $g$, if there is a continuous surjective map $h:M\to N $  such that 
	\begin{eqnarray*}
		h\circ f=g\circ h.
	\end{eqnarray*}
\end{Def}

\begin{Def} 
	Let $X$ be a metric space. Consider a  (piecewise) continuous map $f:X\to X$. If there exists a compact invariant set $\Lambda \subset X $ such that the restriction of $f$ to $\Lambda$ is semi-conjugate to the m-shift map $\sigma $, which is defined as 
	\[\sigma(s)_i=s_{i+1},\]
	then $f$ is said to have an m-type topological horseshoe.  
\end{Def}

For the practical identification of horseshoes within real-world problems, it is necessary to introduce the concept of crossing, as detailed in \cite{Yang2009}. Consider a compact and connected region \( D \subset \mathbb{R}^n \). Let \( B_i \) (\( i = 1, 2, \ldots, m \)) be compact, path-connected subsets in \( D \), each homeomorphic to the unit cube. Denote the boundary of each set \( B_i \) by \( \partial B_i \). Consider a piecewise continuous map \( f: D \rightarrow X \), which is continuous on each compact set \( B_i \).

\begin{Def} \cite{Yang2009}
	For each $B_i, 1\leq i\leq m$, let $B_i^1$ and $B_i^2$ be two fixed disjoint connected nonempty compact subsets (usually pieces of $\partial B_i$) contained in the boundary $\partial B_i$. A connected subset $l$ of $B_i$ is said to be a connection of $B_i^1$ and $B_i^2$ if $l\cap B_i^1\neq \emptyset$ and $l\cap B_i^2\neq \emptyset$. 
\end{Def}

\begin{Def}\label{thd} \cite{Yang2009} 
	Let $l\subset B_i$ be a connection of $B_i^1$ and $B_i^2$. We say that $f(l)$ is crossing $B_j$, if $l$ contains a connected subset $\bar{l}$ such that $f(\bar{l})$ is a connection of $B_j^1$ and $B_j^2$, i.e., $f(\bar{l})\subset B_j$, while $f(\bar{l})\cap B_j^1\neq \emptyset$ and $f(\bar{l})\cap B_j^2\neq \emptyset$. In this case we denote it by $f(l)\mapsto B_j$. Furthermore, if $f(l) \mapsto B_j $ for every connection $l$ of $B_i^1$ and $B_i^2$, then $f(B_i)$ is said to be crossing $B_j$ and denoted by $f(B_i)\mapsto B_j$. To simplify the terminology, we refer to $B_i$ as a crossing block of $B_j$, and for clarity, we call this crossing the dimension one crossing.
\end{Def}

In light of the aforementioned definitions, we recall the following theorem.
\begin{thm}\label{thm1}
	\cite{Yang2004} Suppose that the map $f:D\to \mathbb{R}^n$ satisfies the following assumptions:\\		
	\noindent(1) There exist $m$ mutually path-connected disjoint compact subsets $B_1,B_2,\dots$ and $B_m$ of $D$, the restriction of $f$ to each $B_i$, i.e., $f|_{B_i}$ is continuous.\\		
	\noindent (2) The dimension one crossing relation $f(B_i)\mapsto B_j $ holds for $1\leq i,j\leq m$.\\	
	Then there exists a compact invariant set $K\subset D$, such that $f|_K$ is semi-conjugate to a m-shift map.
\end{thm}

\begin{remark}
	 In the present paper, the crossing structure in Theorem \ref{thm1} is referred to as complete crossing, and the corresponding topological horseshoe can be regarded as a complete topological horseshoe.  We emphasize this point in order to distinguish them from the latter discussion of incomplete corssing and semi-topological horseshoes.
\end{remark}

\section{Subshift of finite type}\label{pre}

When detecting topological horseshoes, situations occasionally arise wherein complete crossing relationships are absent. This necessitates a return to symbolic dynamical systems for further insights. Despite their seeming simplicity, symbolic dynamical systems exhibit a wealth of dynamic behaviors and are fundamental in the study of chaos theory.
In symbolic dynamical systems, subshift of finite type is a powerful tool, and can help us evaluate the topological entropy \cite{Adler1965}, which is introduced by Adler, Konheim, and McAndrew in 1965 to quantify the complexity of symbolic dynamical systems.

A significant tool to describe subshifts of finite type is the non-negative square matrix. For the purposes of this paper, it suffices  to consider $0,1$-square matrix. We will employ the $0,1$-square matrix to define a subshift and subsequently introduce its properties in this section. For more details, readers can refer to \cite{Douglas2021}.

 Firstly, we begin with a succinct definition of a $0,1$-square matrix.
\begin{Def}
	Let $m\geq 2$ be a integer. We say a matrix $A$ is a  $0,1-$squire matrix of order $m$, if 
	\[a_{ij}\in\{0,1\},~i,j=1,2,\dots,m,\]
	where $a_{ij}$ is the element of A in row i and column j.
\end{Def}

Next, we define a subshift determined by a $0,1$-square matrix $A$.

\begin{Def}Let $S_m=\{0,1,\dots,m-1\}$ be the set of nonnegative integers from 0 to $m-1$. Let $\Sigma_m$ be the collection of all one-sided sequences with their elements in $S_m$.	
	Given a $0,1$-square matrix $A$ of order $m$, we can determine a subset $\Sigma_A$ of $\Sigma_m$ as follows:
	\[\Sigma_A=\{s\in \Sigma_m|a_{s_i,s_{i+1}}=1, i\geq 0\}.\] 
	Then, a subshift $\sigma_A:\Sigma_A\to \Sigma_A$ is defined as 
	\[\sigma_A=\sigma|_{\Sigma_A}.\]
	Here, we call $A$ the adjacency matrix of $\sigma_A$.
\end{Def}

In \cite{Douglas2021} and \cite{Zhou1997}, the dynamics of subshifts have been thoroughly investigated. The question of whether a subshift is chaotic hinges on its minimality. To begin with, we revisit the definition of a minimal map.

\begin{Def}\cite{Douglas2021,Zhou1997}
	Let $X$ be a compact metric space, $f:X\to X$ be a continuous map. If
	\[\overline{\mathrm{orb}(x)}=X,~\forall x\in X,\]
	that is to say the orbit of each point in $X$ is dense, then  $f$ is said to be a minimal map.
\end{Def}

A minimal map implies non-chaotic behavior. This leads us to the following theorem.

\begin{thm}\label{lem2}\cite{Douglas2021,Zhou1997}
	Let $k\geq 2$, $A$ be a 0,1-square matrix of order $m$. Suppose $A$ is irreducible, which means for fixed $1\leq i,j\leq m$, there exists $n>0$, such that
	\[(A^n)_{ij}>0,\]
	where $(A^n)_{ij}$ is the element of $A^n$ in row i and column j.
	Then, $\sigma_A$ is chaotic if and only if $\sigma_A$ is not minimal.
\end{thm}

Additionally, there is a theorem characterizing the topological entropy of a subshift.
\begin{thm}\label{thm2}
	\cite{Douglas2021,Zhou1997} 
	Let $h(\sigma_A)$ be the topological entropy of $\sigma_A$, then we have
	\[h(\sigma_A)=\ln\rho(A),\]
	where $A$ is the adjacency matrix of $\sigma_A$, and $\rho(A)$ represents the spectral radius of $A$.
\end{thm}

\section{A generalization of topological horseshoe lemma}\label{ITHT}

In order to generalize the topological horseshoe lemma(Theorem \ref{thm1}), we introduce the definitions of incomplete crossing and semi-topological horseshoe.

\begin{Def}
		Let $D$ be a compact subset of $\mathbb{R}^n$.
		Suppose $f: D \to \mathbb{R}^n$ satisfies the following assumptions:\\
		(1) There exist   $m$ mutually path-connected disjoint compact subsets $B_1,B_2,\dots$ and $B_m$ of $D$, and the restriction of $f$ to each $B_i$, denoted $f|_{B_i}$, is continuous.\\
		(2) For all $i\in\{1,2,\dots,m\}$, there exists $j\in\{1,2,\dots,m\}$, such that $f(B_i)\mapsto B_j$.\\
		(3) There exists a pair $(i,j) \in \{1,2,\dots,m\} \times \{1,2,\dots,m\}$ such that 
		$f(B_i) \cap B_j = \emptyset$, then we say $f$ has a semi-topological horseshoe or semi-horseshoe for simplicity, and such crossing structure is referred to be incomplete crossing.
\end{Def}

To determine whether a semi-topological horseshoe  implies chaos, the following discussion is necessary.

Let $X$ be a metric space, $D$ is a compact subset of $X$, and $f:D\to X$ is a map. 
In the rest part of this section, we make the following assumption:
\begin{ass}\label{ass1}
	There exist $m$ mutually disjoint compact subsets $D_1,D_2,\dots, D_m\subset D$, the restriction of $f$ to each $D_i$, i.e., $f|D_i$ is continuous.
\end{ass}

	Motivated by the adjacency matrix, we introduce the concept of the crossing matrix.
	
	\begin{Def}
		For the sets $D_1,D_2,\dots, D_m\subset D$ and map $f$ mentioned above, we define a $0,1$-matrix $A=(a_{ij})$ as follows:
		\begin{eqnarray*}
			a_{ij}=\left\{\begin{array}{l}
				1,~\mathrm{if}~f(B_i)\mapsto B_j,\\
				0,~\mathrm{otherwise}.
			\end{array}\right.
		\end{eqnarray*}
		$A$ is said to be the crossing matrix of $f$.
	\end{Def}

	For the sake of discussion, we make another assumption.
	
	\begin{ass}\label{ass2}
		The square matrix $A$ mentioned in this section is irreducible, that is to say for any fixed $1\leq i,j\leq m$, there exists $n>0$, such that
		\[(A^n)_{ij}>0,\]
		where $(A^n)_{ij}$ is the element of $A^n$ in row i and column j.
	\end{ass}

\begin{remark}
	If the condition in Assumption \ref{ass2} fails, for instance there exist $0\leq i,j \leq m$, such that $(A^n)_{ij}=0, \forall n\geq 0$. That is to say, $f^n(B_i)\cap B_j\equiv \emptyset,~\forall n\geq 0$, and the intersection of $B_i$ and the  chaotic invariant set is an empty set. So, Assumption \ref{ass2} is necessary.
\end{remark}

\begin{Def}
	Let $\gamma$ be a compact set contained in $D$. If $\gamma\cap D_i\neq \emptyset$, then we set $\gamma_i=\gamma\cap D_i$. Define $F_A$ as follows:
	\[F_A=\{\gamma\subset D |\gamma ~is ~comact,~f(\gamma_i)\cap D_j\neq \emptyset, \forall a_{ij}=1\}.\]
	If $\gamma\subset F_A$, then $f(\gamma_i)=f(\gamma\cap D_i)$ is a compact set contained in $D$. For all $(i,j): a_{ij}=1$, we have
	\[\emptyset\neq f(\gamma_i)\cap D_j=:\gamma_{ij}\subset D_j.\]
	Furthermore,
	\[f(\gamma_{ij})\cap D_k\neq\emptyset,~\forall a_{jk}=1,\]
	which means $f(\gamma_i)\in F_A.$
	
	We call $F_A$ a $f_A$-connected family with respect to $D_1,D_2,\dots$ and  $D_m$. 
\end{Def}

\begin{Lem}\label{lem1}
Suppose that there is a $f_A$-connected family $F_A$ with respect to $D_1,D_2,\dots$ and $D_m$.
Then, for every $s\in\Sigma_A$,
	\[s=\{s_0,s_1,\dots,s_n,\dots\}, s_i\in S_m,\]
	there exists a point $x\in D$, such that $f^n(x)\in D_{s_n }$.
\end{Lem}

\noindent{\bf Proof.}
	Fix $\gamma\in F_A$, set $\gamma_{s_0}=\gamma\cap D_{s_0}$.
	Since $s\in \Sigma_A$, we have $a_{s_0s_1}=1$, it follows that
	\[f(\gamma_{s_0})\in F_A.\]
	Set $\gamma_{s_0s_1}=f(\gamma_{s_0})\cap D_{s_1}\neq \emptyset$, then \[f(\gamma_{s_0s_1})\in F_A.\]
	Inductively, for $n\geq 0$, we have
	\[\gamma_{s_0s_1\dots s_n}=f(\gamma_{s_0s_1\dots s_{n-1}})\cap D_{s_n}\neq \emptyset,\]
	hence
	\[f(\gamma_{s_0s_1\dots s_n})\in F_A.\]
	For $n\geq 0$, set
	\[\Gamma_n^1=\gamma_{s_0s_1\dots s_n}.\]	
	Set inductively for $j\geq 1$:
	\[\Gamma_n^{j+1}=\{x\in\Gamma_n^j|f(x)\in\Gamma_{n+1}^j\}.\]	
	It is obvious that
	\[\Gamma_n^{j+1}\subset \Gamma_n^j \subset \Gamma_n^{j-1},\]
	hence it is evident that $\{\Gamma_n^j\}_{j=1}^\infty$ forms a nested sequence of compact sets, leading to the definition:
	\[\Gamma(n):=\bigcap\limits_{j=1}^\infty \Gamma_n^j.\]
	In the following, we  aim to prove  $f(\Gamma_n^{j+1})=\Gamma_{n+1}^j$ by induction.
	
	Firstly, from
	\[\Gamma_{n+1}^1=\gamma_{s_0\dots s_{n+1}}=f(\gamma_{s_0\dots s_n})\cap D_{s_{n+1}},\]
	\[\Gamma_n^2 =\{x\in\Gamma_n^1|f(x)\in \Gamma_{n+1}^1\}=\{x\in\gamma_{s_0\dots s_n}|f(x)\in D_{s_{n+1}}\},\]
	we have $f(\Gamma_n^2)=\Gamma_{n+1}^1$.
	
	Suppose that $f(\Gamma_n^{j+1})=\Gamma_{n+1}^j$. We will prove $f(\Gamma_n^{j+2})=\Gamma_{n+1}^{j+1}$ next.
	\[\Gamma_{n+1}^{j+1}=\{x\in \Gamma_{n+1}^j|f(x)\in\Gamma_{n+2}^j\}\]
	\[\Gamma_n^{j+2}=\{x\in \Gamma_n^{j+1}|f(x)\in\Gamma_{n+1}^{j+1}\}\]
	\[f(\Gamma_n^{j+2})=f(\Gamma_n^{j+1})\cap \Gamma_{n+1}^{j+1}=\Gamma_{n+1}^j \cap\Gamma_{n+1}^{j+1}=\Gamma_{n+1}^{j+1}. \]
	The last equation is because $\Gamma_{n+1}^{j+1}\subset\Gamma_{n+1}^j$.
	Thus,
	\[f(\Gamma(n))=f(\bigcap\limits_{j=1}^\infty \Gamma_n^j)=f(\bigcap\limits_{j=2}^\infty \Gamma_n^j)=\bigcap\limits_{j=2}^\infty f(\Gamma_n^j)=\bigcap\limits_{j=2}^\infty \Gamma_{n+1}^{j-1}=\bigcap\limits_{j=1}^\infty \Gamma_{n+1}^j=\Gamma(n+1).\]
	It follows that for all $x\in\Gamma(0)$, $f^n(x)$ is well defined, and 
	\[f^n(x)\in \gamma_{s_0s_1\dots s_n}\subset D_{s_n}.\]
\begin{flushright}
	\opensquare
\end{flushright}

\begin{thm}\label{ITH}
	Suppose that there exists a $f_A$-connected family $F_A$ with respect to $D_1,\dots,D_m$. Then there exists a compact invariant set $K\subset D$, such that $f|K$ is semi-conjugate to a subshift $\sigma_A$, where $A$ is the corresponding adjacency matrix.
\end{thm}

\noindent{\bf Proof.}
	The proof given here follows the idea from \cite{Yang2004}.
	
	Firstly, we will demonstrate the existence of a point $\bar{x}\in D$ such that for every $s=\{s_0,s_1,\dots,s_n,\dots\}\in \Sigma_A$, there exists a point $x\in\omega(\bar{x})$ satisfying $f^n(x)\in D_{s_n}$.
	
	To achieve this objective, we select a sequence $\bar{s}\in \Sigma_A$, which is the initial condition of a dense orbit on $\Sigma_A$. Given the irreducibility of 
	$A$, from Theorem \ref{lem2}, we know that $\sigma_A$ is chaotic if and only if $\sigma_A$ is not minimal. When $\sigma_A$ is minimal, for every point $s\in\Sigma_A$, $\mathrm{orb}(s)$ is dense in $\Sigma_A$. Otherwise, if $\sigma_A$ is chaotic, then there exists a orbit dense in $\Sigma_A$. Therefore, the sequence $\bar{s}$ mentioned earlier indeed exists.
	To the sequence
	\[\bar{s}=\{\bar{s_0},\bar{s_1},\dots,\bar{s_n},\dots\},\]
	we associate the corresponding sequence of symbolic sets
	\[D_{\bar{s}}=\{D_{\bar{s_0}},D_{\bar{s_1}},\dots,D_{\bar{s_n}},\dots\}.\]
	By Lemma \ref{lem1}, there exists $\bar{x}\in D_{\bar{s}}$ such that $f^i(\bar{x})\in D_{\bar{s_i}},~\forall i\geq 0$.
	It is easy to see there is a sequence $\{\varepsilon_k\}$ satisfying $\lim\limits_{k\to\infty}\varepsilon_k= 0$, such that for all $s\in \Sigma_A$ satisfying
	\[d(s,\bar{s})=\sum\limits_{i=0}^\infty \frac{1}{2^i} \frac{|s_i-\bar{s_i}|}{1+|s_i-\bar{s_i}|} <\varepsilon_k,\] 
	we have $ s_i=\bar{s_i},~ 0\leq i \leq k$.
	Consider $s\in\Sigma_A$. Since $\mathrm{orb}(\bar{s})$ is dense in $\Sigma_A$, for each $\varepsilon_k$ mentioned above, there exists a number $n_k>0$, such that 
	\[d(s,\sigma^{n_k}_A(\bar{s}))< \varepsilon_k,\] 
	which implies $s_i=(\sigma^{n_k}_A(\bar{s}))_i,~0\leq i\leq k.$
	Therefore, \[f^{i+n_k}(\bar{x})\in D_{s_i},~ 0\leq i\leq k. \]
	Now since  $\{f^{0+n_k}(\bar{x})\}_{n_k>0}\subset D_{s_0}$, there exists a subsequence of $\{n_k\}$, say $\{n_k^0\}$, such that 
	\[f^{0+n_k^0}(\bar{x})\to x_0\in D_{s_0}~ (n_k^0\to\infty ).\]
	Similarly, for $\{f^{1+n_k^0}(\bar{x})\}_{n_k^0>0}\subset D_{s_1}$, there exists a subsequence of $\{n_k^0\}$, say $\{n_k^1\}$, such that 
	\[f^{1+n_k^1}(\bar{x})\to x_1\in D_{s_1}~ (n_k^1\to\infty ).\]
	Following this process recursively, there is a subsequence $\{n_k^i\}\subset \{n_k^{i-1}\}(\forall i\geq0)$, such that 
	\[f^{i+n_k^i}(\bar{x})\to x_i\in D_{s_i}, ~(n_k^i\to \infty).\]
	It follows that for each $i>0$:
	\begin{eqnarray*}
		\lim\limits_{n_k^i\to\infty} f^{i+n_k^i} (\bar{x}) &= \lim\limits_{n_k^i\to\infty} f^i(f^{n_k^i}(\bar{x}))\\
		&=\lim\limits_{n_k^0\to\infty} f^i(f^{n_k^0}(\bar{x}))\\
		&=f^i(\lim\limits_{n_k^0\to\infty} f^{n_k^0}(\bar{x}))\\
		&=f^i(x_0),
	\end{eqnarray*}
	which implies that $f^i(x_0)\in D_{s_i}$. Since $f^{0+n_k^0}(\bar{x})\to x_0$, we have $x_0\in\omega(\bar{x})$. In particular,  setting  $s=\bar{s}$, there is a point $x^0\in\omega(\bar{x})$, such that $f^i (x^0)\in D_{\bar{s_i}}$.
	
	Secondly, we will establish the existence of $\bar{x}\in\bigcup\limits_{i=1}^m D_i\subset D$ such that \( f^i(\bar{x}) \in D_{\bar{s_i}} \) for \( i \geq 0 \), and moreover, \(\bar{x} \in \omega(\bar{x})\).
	
	Consider a collection of compact sets defined as follows:
	\[\{\omega(x)|x\in\bigcup\limits_{i=1}^m D_i\subset D: f^i(x)\in D_{\bar{s_i}},~i\geq 0\},\]
	which is partially ordered by inclusion.
	According to Hausdorff's Maximal Theorem, there exists a maximal totally ordered subcollection $B$.
	Let $M=\bigcap\limits_{\Lambda\in B} \Lambda$, then  $M$ is a compact invariant set.
	For all $\bar{x}\in M$, we have $\bar{x}\in \bigcup\limits_{i=1}^m D_i, f^i(\bar{x})\in D_{\bar{s_i}},~\forall i\geq 0$ and $\bar{x}\in \omega(\bar{x})$.
	
	Thirdly, we will demonstrate that for each \( x \in \omega(\bar{x}) \), the orbit \( \{f^i(x) \mid i \geq 0\} \) is contained within \( \bigcup_{i=1}^m D_i \subset D \).
	
	Suppose not, then there exists some $j>0$, such that 
	\[f^j(x)\notin\bigcup\limits_{i=1}^m D_i.\]
	Since $\bigcup\limits_{i=1}^m D_i$ is closed, there must exist an open neighborhood $N_j$ of $f^j(x)$ satisfying $N_j\cap \bigcup\limits_{i=1}^m D_i=\emptyset$.
	Given  $\{f^i(x)|i\geq 0\}\subset\bigcup\limits_{i=1}^m D_i$, it follows that $\{f^i(x)|i\geq 0\}\cap N_j=\emptyset$, where $f^j(x)\in N_j$, implying  $f^j(x)\notin \omega(\bar{x})$. However, since $x\in \omega(\bar{x})$, we have $f^j(x)\in\omega( \bar{x})$, leading to a contradiction.
	
	Fourthly, define a map $g:\omega(\bar{x})\to \Sigma_A$ as follows:
	\[g(x)=s=\{s_0,s_1,\dots,s_n,\dots\}, ~\mathrm{if}~ f^i(x)\in D_{s_i}, \forall i\geq 0.\]
    Due to the continuity of $f$ on each $D_i$, it is evident  that $g$ is continuous.
	
	Finally, it is straightforward to verify that
	\[g\circ f(x)=\sigma_A\circ g(x),~ \forall x\in \omega(\bar{x}).\]
	The proof is complected.	
\begin{flushright}
	\opensquare
\end{flushright}

\begin{col}\label{col1}
	In Theorem \ref{ITH}, if $\sigma_A$ is not minimal, then $f$ is chaotic, and furthermore we have
	\[h(f)\geq \ln(\rho(A)).\]
\end{col}

\begin{remark}

	We additionally present a sufficient condition for the detection of semi-horseshoes and provide an rough estimation of the topological entropy.
\end{remark}	

\begin{Prop}
		Suppose the crossing matrix of $f$, noted as $A$, is irreducible, and there are $k$ rows in $A$, whose elements are all 1:
		\[a_{J_ji}=1, ~1\leq j\leq k, ~1\leq i\leq m.\]
		Define $N_{J_j},~1\leq j\leq k$ as follows:
		\[N_{J_j}=\#\{a_{iJ_j}=1,~1\leq i\leq m \}~1\leq j\leq k,\]
		where $\#\{\cdot\}$ represents the card of the set. 
		Then there exist $\sum\limits_{j=1}^k N_{J_j}$ crossing blocks with respect to $f^2$, which means
		\[h(f)\geq \frac12 \ln(\sum\limits_{j=1}^k N_{J_j}).\]
\end{Prop}

\section{Application of semi-topological horseshoe theory in dynamical systems}\label{Examples}

As an application of the  theory above, we will discuss the perturbed Duffing system and a polynomial system proposed by Chen in \cite{Chen1999} in this section. To begin with, let us focus on the simplest cases -- topological horseshoes involving two blocks.

\subsection{Horseshoes involving two blocks}\label{TwoBlocks}

\begin{figure}[h]
	\centering
	\includegraphics[height=8.5cm]{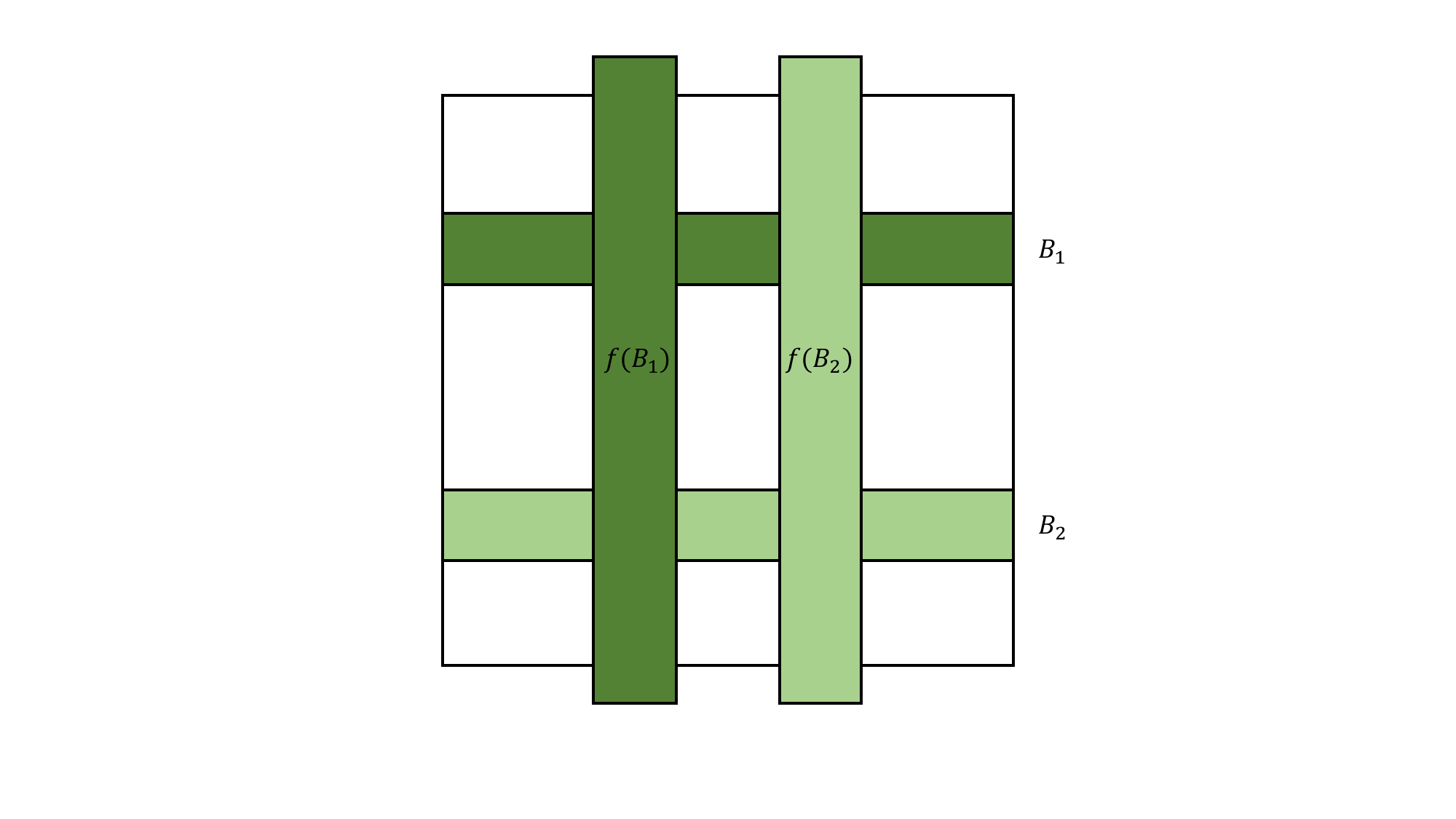}
	\caption{The crossing matrix corresponding to this case  is
		$\left(
			\begin{array}{cc}
			1 & 1\\
			1 & 1
			\end{array}
			\right)$.}
	\label{ppt-1}
\end{figure}

The complete topological horseshoe comprising two crossing blocks (Figure \ref{ppt-1}) is commonly used to verify chaos, so we won't go into details.  It is well-known that the map \( f \) in this scenario is semi-conjugate to a 2-shift map, implying that
\[ h(f) \geq \ln 2, \]
where \( h(f) \) denotes the topological entropy of \( f \).

\begin{figure}[h]
	\centering
	\includegraphics[height=8.5cm]{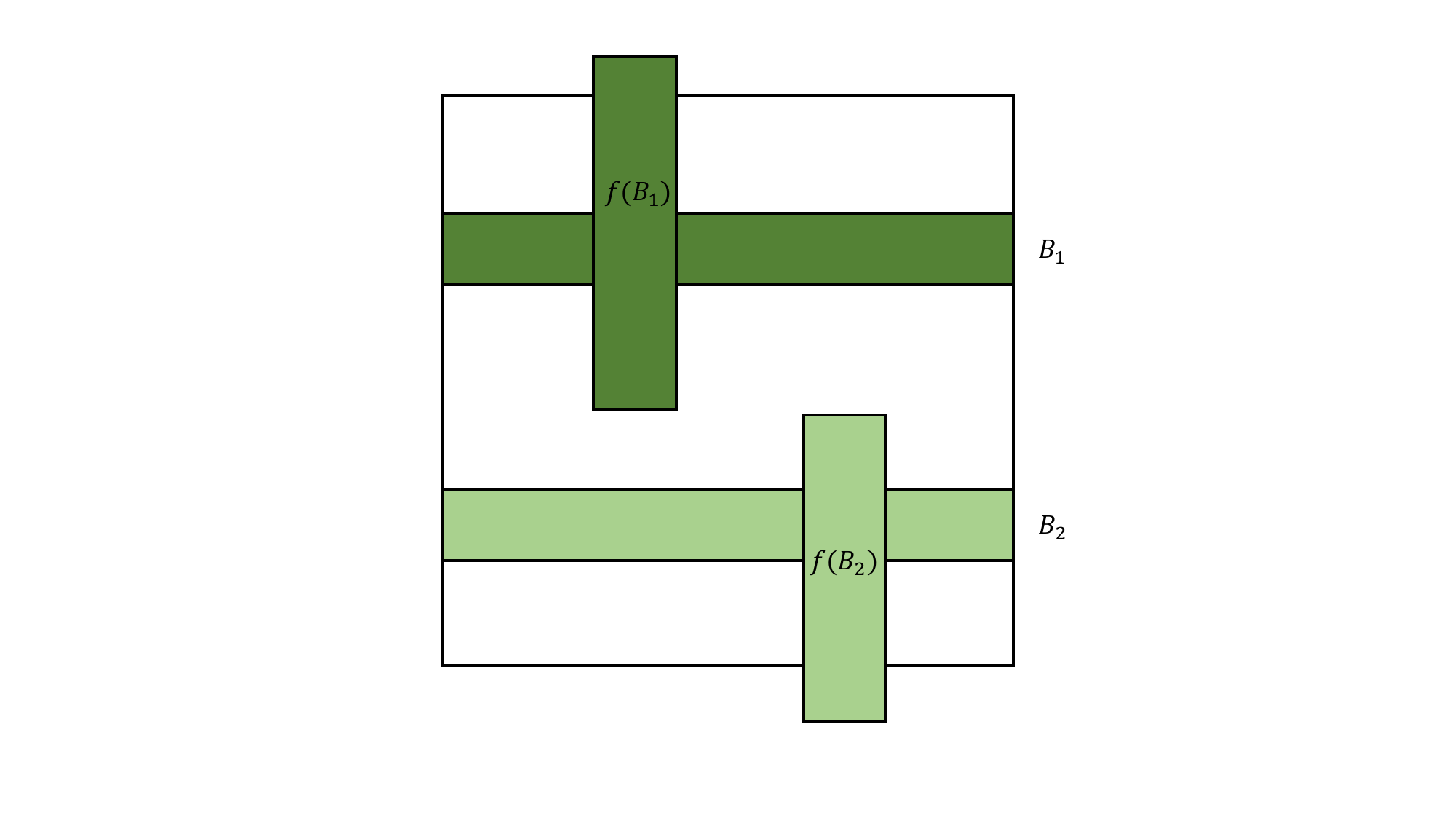}
	\caption{The crossing matrix $A_1$ corresponding to this case  is
		$\left(
		\begin{array}{cc}
		1 & 0\\
		0 & 1
		\end{array}
		\right)$.}		
	\label{ppt-7}
\end{figure}
\begin{figure}[h]
	\centering
	\includegraphics[height=8.5cm]{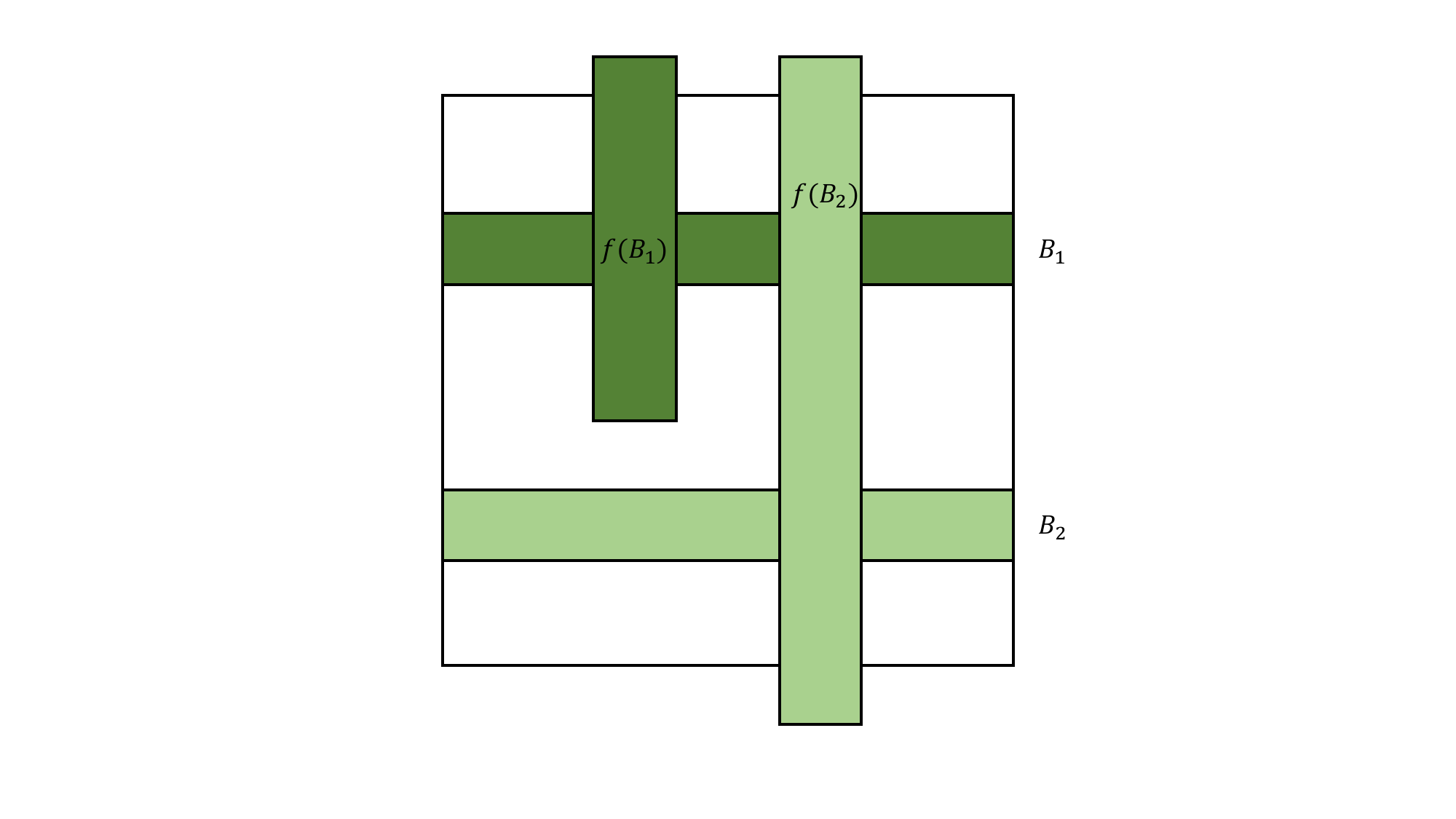}
	\caption{The crossing matrix $A_2$ corresponding to this case  is 
		$\left(
		\begin{array}{cc}
			1 & 0\\
			1 & 1
		\end{array}
		\right)$.}		
		\label{ppt-6}
\end{figure}

In Figure \ref{ppt-7} and Figure \ref{ppt-6}, both the crossing matrices $A_1$ and $A_2$ are reducible, so these two kinds of semi-topological horseshoe will not lead to chaos.

\begin{figure}[h]
	\centering
	\includegraphics[height=8.5cm]{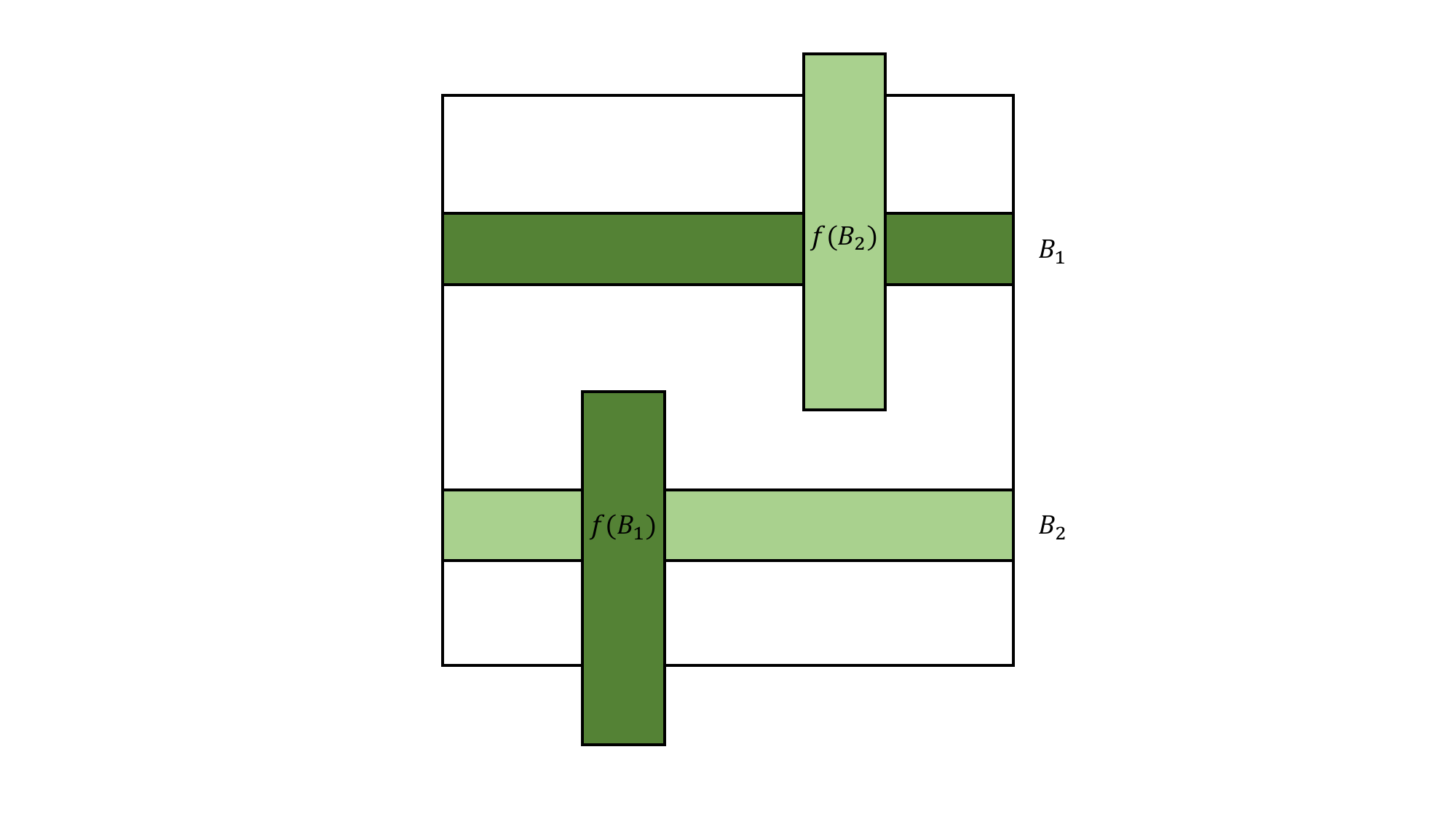}
	\caption{The crossing matrix $A_3$ corresponding to this case  is
			$\left(
		\begin{array}{cc}
			0 & 1\\
		1 & 0
		\end{array}
		\right)$.}		
		\label{ppt-8}
\end{figure}

This time, $A_3$ is irreducible, but $\sigma_{A_3}$ is minimal, since there are only two points in $S|_{A_3}$:
\[(0,1,0,1,\dots)~\&~(1,0,1,0,\dots).\]

\begin{figure}[h]
	\centering
	\includegraphics[height=8.5cm]{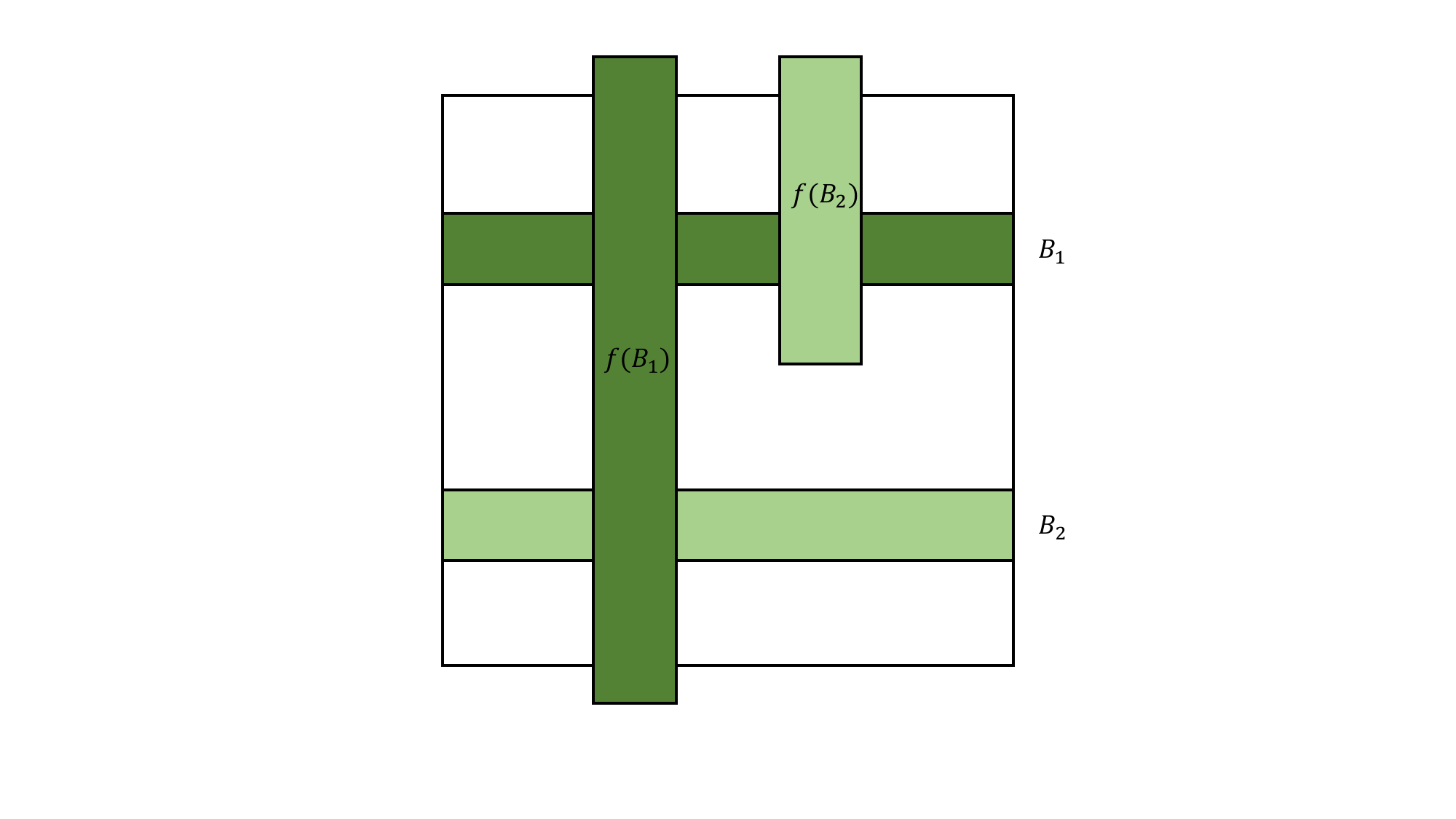}
	\caption{The crossing matrix corresponding to this case is 
		$\left(
		\begin{array}{cc}
			1 & 1\\
			1 & 0
		\end{array}
		\right)$.}		
		\label{ppt-2}
\end{figure}

The next case (Figure \ref{ppt-2}) presents an interesting scenario.

From Corollary \ref{col1}, we have the corollary below.

\begin{col}\label{col2}
	If the crossing matrix of $f$ is \[A=\left(
	\begin{array}{cc}
		1 & 1\\
		1 & 0
	\end{array}
	\right),\]
	then 
	\[h(f)\geq \ln \frac{1+\sqrt{5}}{2}.\]
\end{col}

\noindent{\bf Proof.}
	he crossing matrix of $f$ is \[A=\left(
	\begin{array}{cc}
		1 & 1\\
		1 & 0
	\end{array}
	\right),\] whose eigenvalues, can be easily obtained:
	\[\lambda_1=\frac{1-\sqrt{5}}{2}, ~\lambda_2=\frac{1+\sqrt{5}}{2},\]
	and the spectral radius of $A$ is $\frac{1+\sqrt{5}}{2}$. So, $h(\sigma_A)=\frac{1+\sqrt{5}}{2}$, which implies that
	\[h(f)\geq \ln \frac{1+\sqrt{5}}{2}.\]
\begin{flushright}
	\opensquare
\end{flushright}

\begin{remark}\label{remark1}
	
	For the case depicted in Figure \ref{ppt-2}, we can also study this type of semi-horseshoe through map iteration. For convenience, we introduce the following notations (Figure \ref{ppt-3}):
	\[B_{11}:=f(B_1)\cap B_1,~ B_{12}:=f(B_1)\cap B_2, ~B_{21}:=f(B_2)\cap B_1.\]
	
	\begin{figure}[h]
		\centering
		\includegraphics[height=8.5cm]{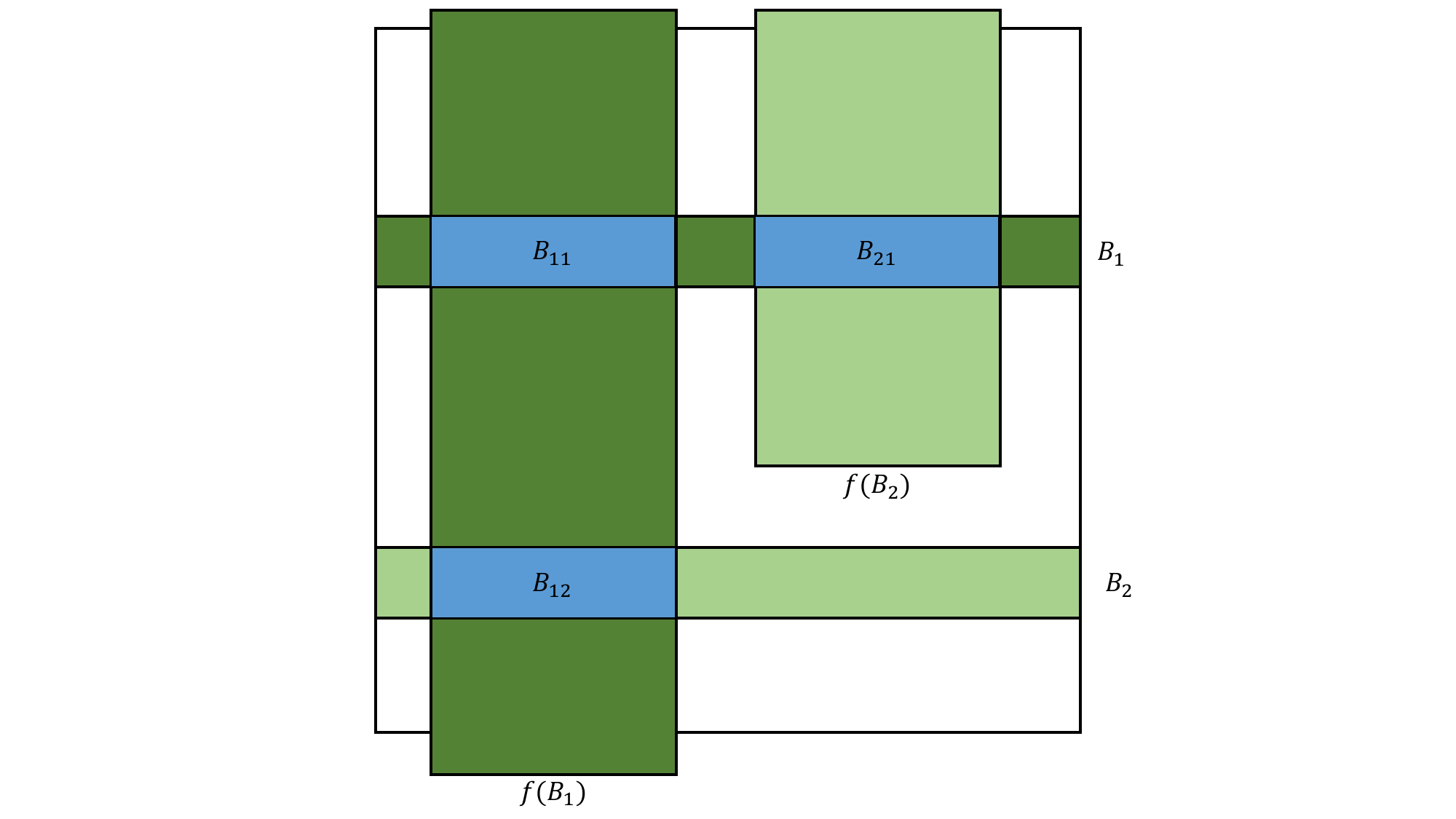}
		\caption{Some notations.}
		\label{ppt-3}
	\end{figure}
	To get the images $f^2(B_1)$ and $f^2(B_2)$, we only have to consider $f(B_{11}),f(B_{12})$ and $f(B_{21})$. Together with the preimages of $B_{11}, B_{12}$ and $B_{21}$, we find two crossing blocks $f^{-1}(B_{11})$ and $f^{-1}(B_{21})$ with respect to $f^2$ (Figure \ref{ppt-4}), which implies that 
	\[h(f^2)\geq \ln2,\] 
	and it follows
	\[h(f)\geq \frac12\ln2.\]
	\begin{figure}[h]
		\centering
		\includegraphics[height=10cm]{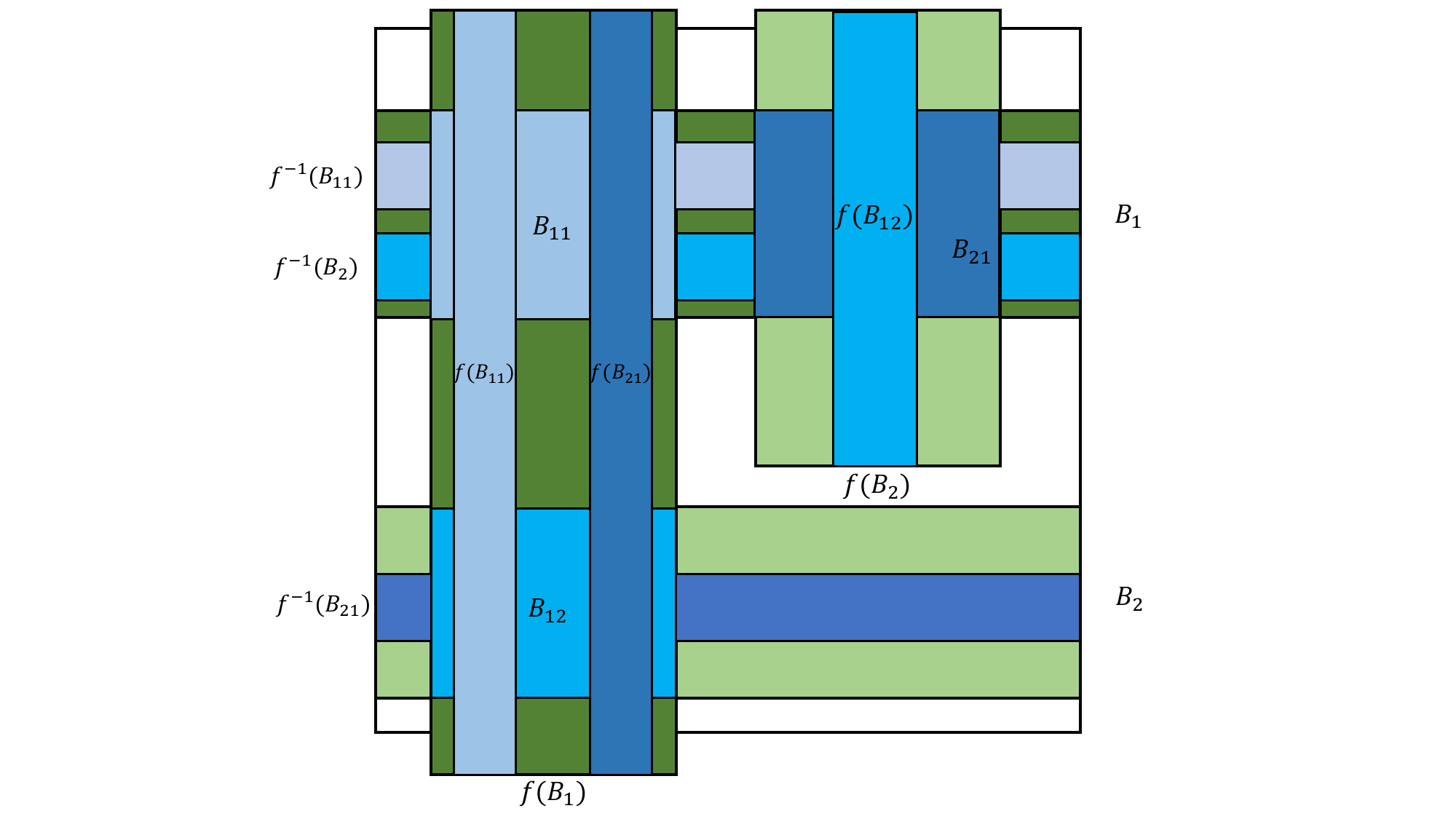}
		\caption{Two crossing blocks with respect to $f^2$.}
		\label{ppt-4}
	\end{figure}
	However, this is not the end. When increasing the iterations of $f$, we observe that the lower bound of the topological entropy of $f$ also increases. When considering $f$, we can find two blocks, whose images under $f$ cross $B_1$, and one of their image crosses $B_2$.
	This implies the presence of two crossing blocks when considering $f^2$.
	Similarly, for considering $f^2$, we can find three blocks, whose images under $f^2$ cross $B_1$, and two of their image crosses $B_2$,
	which implies that we will find 3 crossing blocks when considering $f^2$.
	By simple induction, if we find $x_n$ blocks, whose images under $f^n$ cross $B_1$, and $y_n$ of their image crosses $B_2$ when considering $f^n$,
	then, for \( f^{n+1} \), we will find $x_n+y_n$ blocks, whose images under $f^{n+1}$ cross $B_1$, and $x_n$ of their image crosses $B_2$. This implies the existence of \( x_n \) crossing blocks with respect to \( f^{n+1} \).
	In other words, we obtain the following sequence defined by the recurrence relation:

	\begin{eqnarray*}
		\left\{\begin{array}{l}
			x_{n+1}=x_n+y_n,\\
			y_{n+1}=x_n,\\
			x_1=2,\\
			y_1=1.
		\end{array}\right.
	\end{eqnarray*}	
	Focusing on $\{x_n\}$, we get a Fibonacci sequence:
	\begin{eqnarray*}
		\left\{\begin{array}{l}
			x_{n+2}=x_n+x_{n+1},\\
			x_1=2,\\
			x_2=3
		\end{array}\right.
	\end{eqnarray*}	
	From the general term formula of the Fibonacci sequence, we have 
	\[x_n=\frac{1}{\sqrt{5}} [(\frac{1+\sqrt{5}}{2})^{n+2}-(\frac{1-\sqrt{5}}{2})^{n+2}].\]
	From the discussion above, we deduce  that there exist $x_{n-1}$ crossing blocks with respect to $f^n$, so we have
	\[h(f)\geq\frac{1}{n}\ln\{\frac{1}{\sqrt{5}} [(\frac{1+\sqrt{5}}{2})^{n+1}-(\frac{1-\sqrt{5}}{2})^{n+1}]\}\]
	The term on the right is monotonically increasing with an upper bound, which can be obtained through simple analysis:
	\begin{eqnarray*}
		\lim\limits_{n\to\infty}\frac{1}{n}\ln\{\frac{1}{\sqrt{5}} [(\frac{1+\sqrt{5}}{2})^{n+1}-(\frac{1-\sqrt{5}}{2})^{n+1}]\}=\ln\frac{1+\sqrt{5}}{2}.
	\end{eqnarray*}
	Thus, we can deduce that \[h(f)\geq\ln\frac{1+\sqrt{5}}{2}.\]
\end{remark}

\subsection{Horseshoe and semi-horseshoe in the perturbed Duffing system}\label{Duffing}

The case illustrated in Figure \ref{ppt-2} is not merely an isolated example. We encounter the aforementioned crossing structure when studying the Duffing system with periodic perturbation, which is revisited as follows:
\begin{eqnarray*}\label{Duffing-2}
	\ddot{x} + \delta \dot{x} - x +  x^3 = \gamma \cos(\omega t).
\end{eqnarray*} 
Introducing \( y = \dot{x} \), the system can be reformulated into a first-order form:
\begin{eqnarray*}\label{Non-autonomous}
	\left\{\begin{array}{l}
		\dot{x}=y,\\
		\dot{y}=x-x^3-\delta y+\gamma \cos\omega\theta.
	\end{array}\right.
\end{eqnarray*}
To further analyze the system, we introduce a phase variable \( \theta \), yielding the autonomous system:
\begin{eqnarray}\label{Duffing3}
	\left\{\begin{array}{l}
		\dot{x}=y,\\
		\dot{y}=x-x^3-\delta y+\gamma \cos\omega\theta,\qquad(x,y,\theta)\in\mathbb{R}^2\times S^1,\\
		\dot{\theta}=1.
	\end{array}\right.
\end{eqnarray}

Here, we focus on the case with the parameters set in \cite{Guckenheimer1984}:
\[\omega=1,\delta=0.25.\]

In \cite{Guckenheimer1984}, the Melnikov method \cite{Melnikov1963} along with the Smale-Birkhoff theorem \cite{Smale1965,Smale1967} affirm that the perturbed Duffing system exhibits chaotic dynamics under sufficiently small external forces. However, the precise threshold defining "sufficiently small" remains a current area of investigation.
In this subsection, we examine two specific scenarios where $\gamma$ takes on distinct magnitudes.

\begin{figure}[h]
	\centering
	\includegraphics[height=7cm]{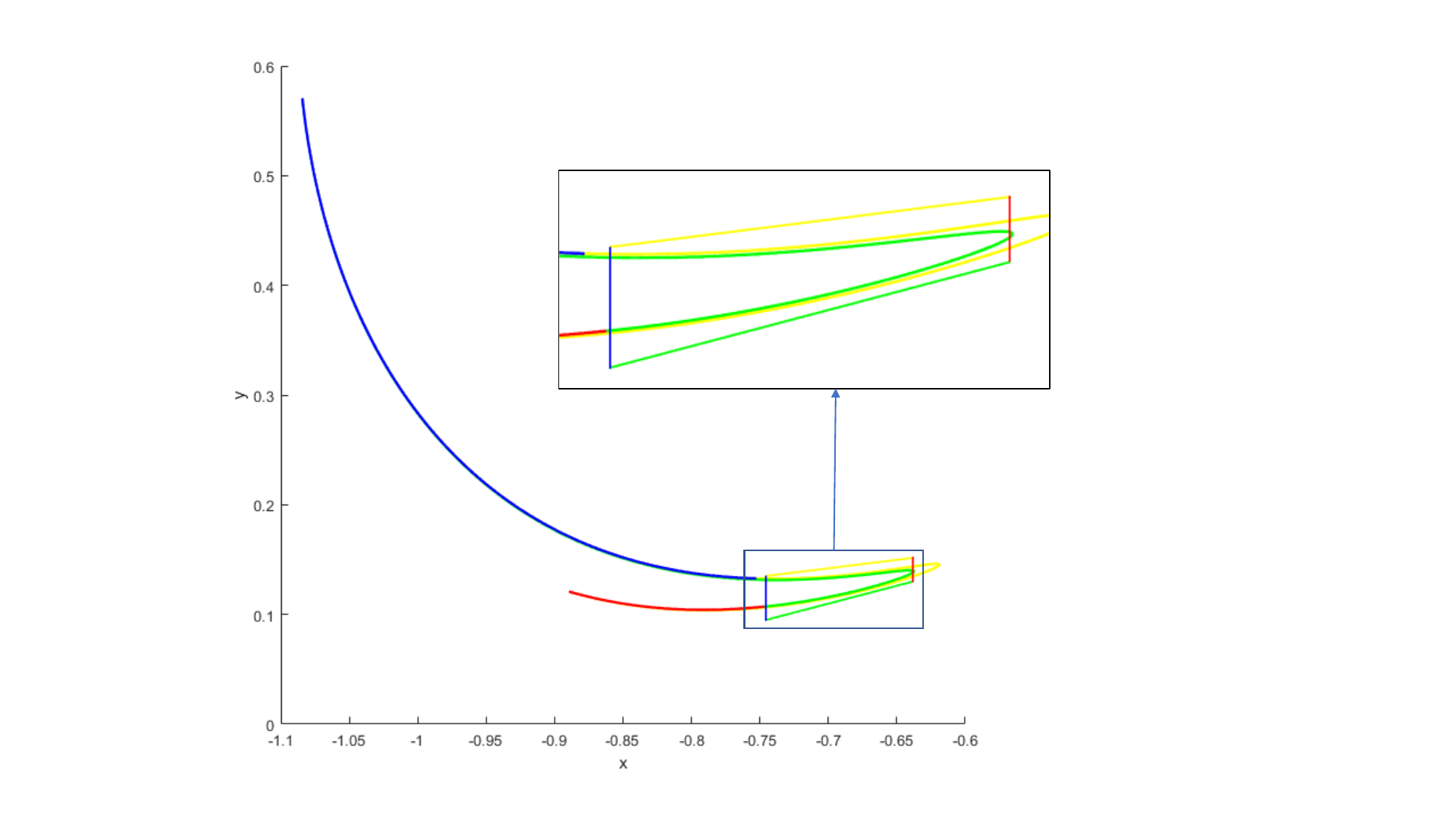}
	\caption{Smale horseshoe ($\gamma$=1.183).}
	\label{1.183}
\end{figure}

When $\gamma=1.183$, Figure \ref{1.183}  shows that the Smale horseshoe exists, indicating that the Poincaré map $P_{\gamma=1.183}$ is semi-conjugate to a 2-shift map. Consequently,
\[h(P_{\gamma=1.183})\geq \ln 2.\]

\begin{figure}[h]
	\centering
	\includegraphics[height=10cm]{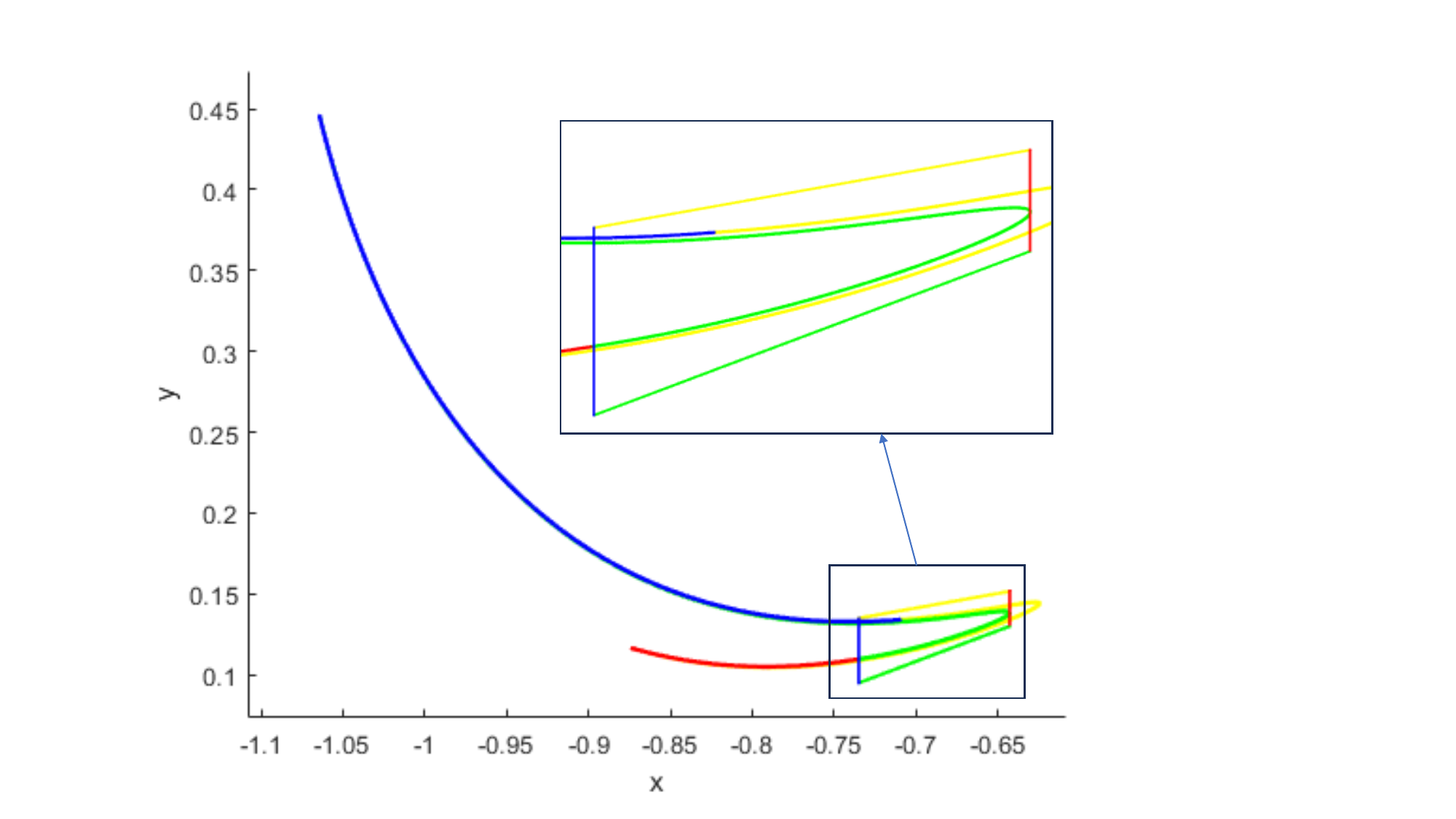}
	\caption{Smale horseshoe disappears($\gamma$=1.184).}
	\label{1.184}
\end{figure}

\begin{figure}[h]
	\centering
	\includegraphics[height=10cm]{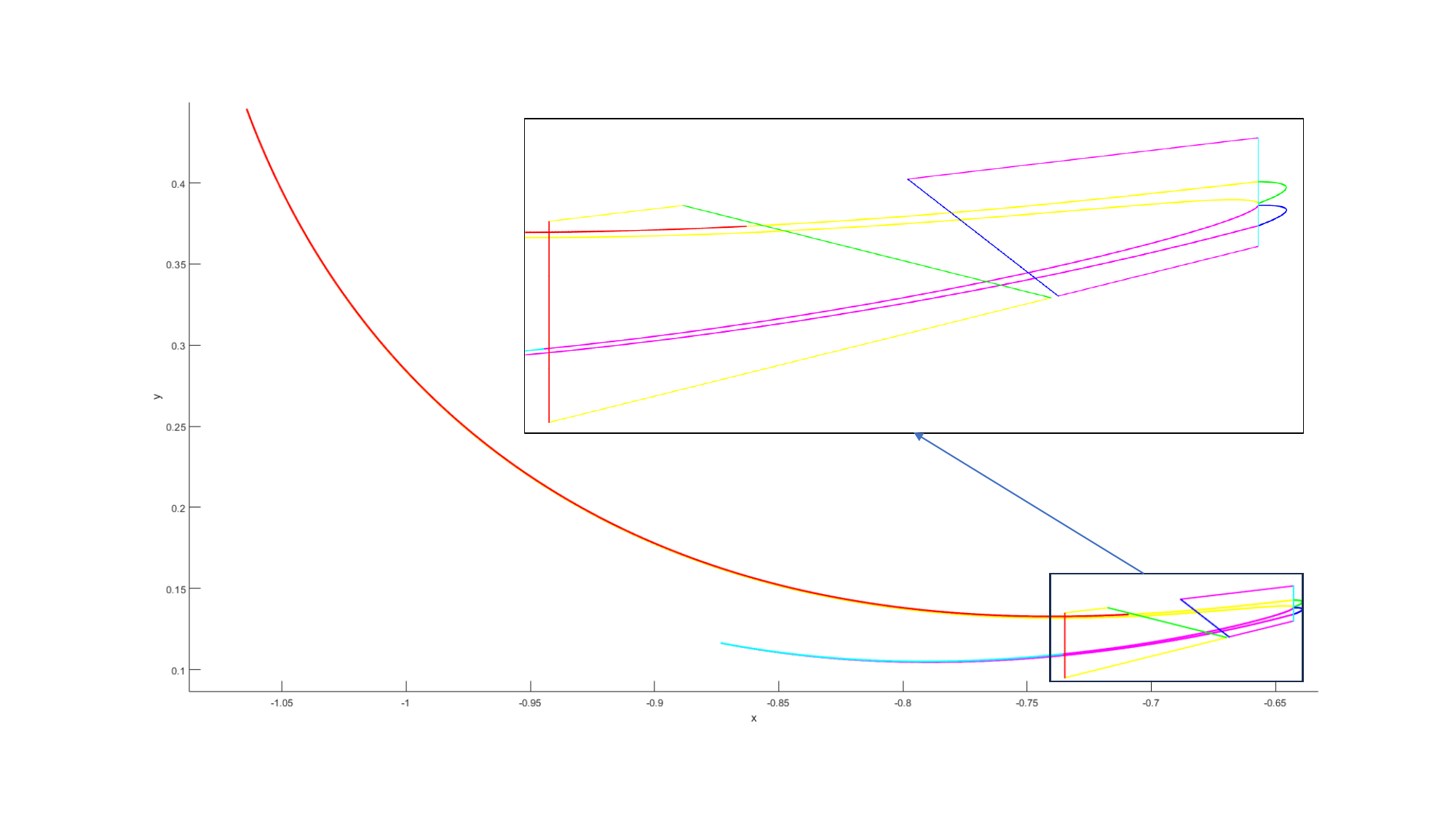}
	\caption{Semi-horseshoe ($\gamma$=1.184).}
	\label{1.184-3}
\end{figure}

When $\gamma$ increases to $1.184$, the Smale horseshoe no longer exists (Figure \ref{1.184}), and directly identifying the topological horseshoe becomes challenging. In the subsequent figures, we employ the Runge-Kutta method to visualize the mapping of blocks under the Poincaré map, where each side of the quadrilaterals is mapped to the side of the corresponding color.
In Figure \ref{1.184-3}, when $\gamma=1.184$, we find a semi-horseshoe.

\begin{figure}[h]
	\centering
	\includegraphics[height=8.5cm]{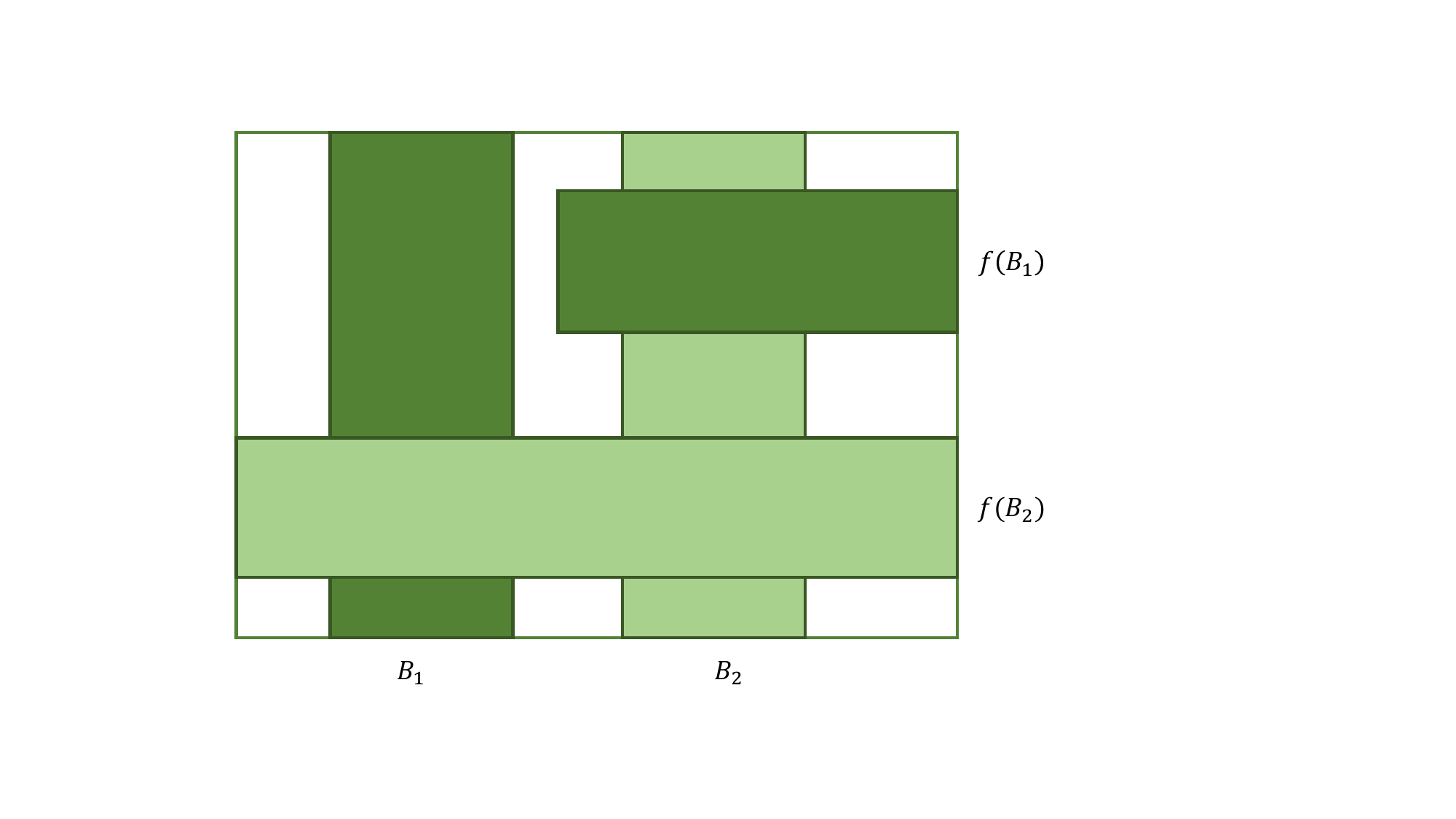}
	\caption{Illustration of incomplete crossing when $\gamma$=1.184.}
	\label{}
\end{figure}

From Corollary \ref{col2}, we can deduce that the topological entropy of the first return map $P$ of the perturbed Duffing system when $\gamma=1.184$  is no less than $\ln(\frac{1+\sqrt{5}}{2})$.

\subsection{Semi-topological horseshoe in a polynomial system}\label{Chen}

\begin{figure}[h]
	\centering
	\includegraphics[height=8.5cm]{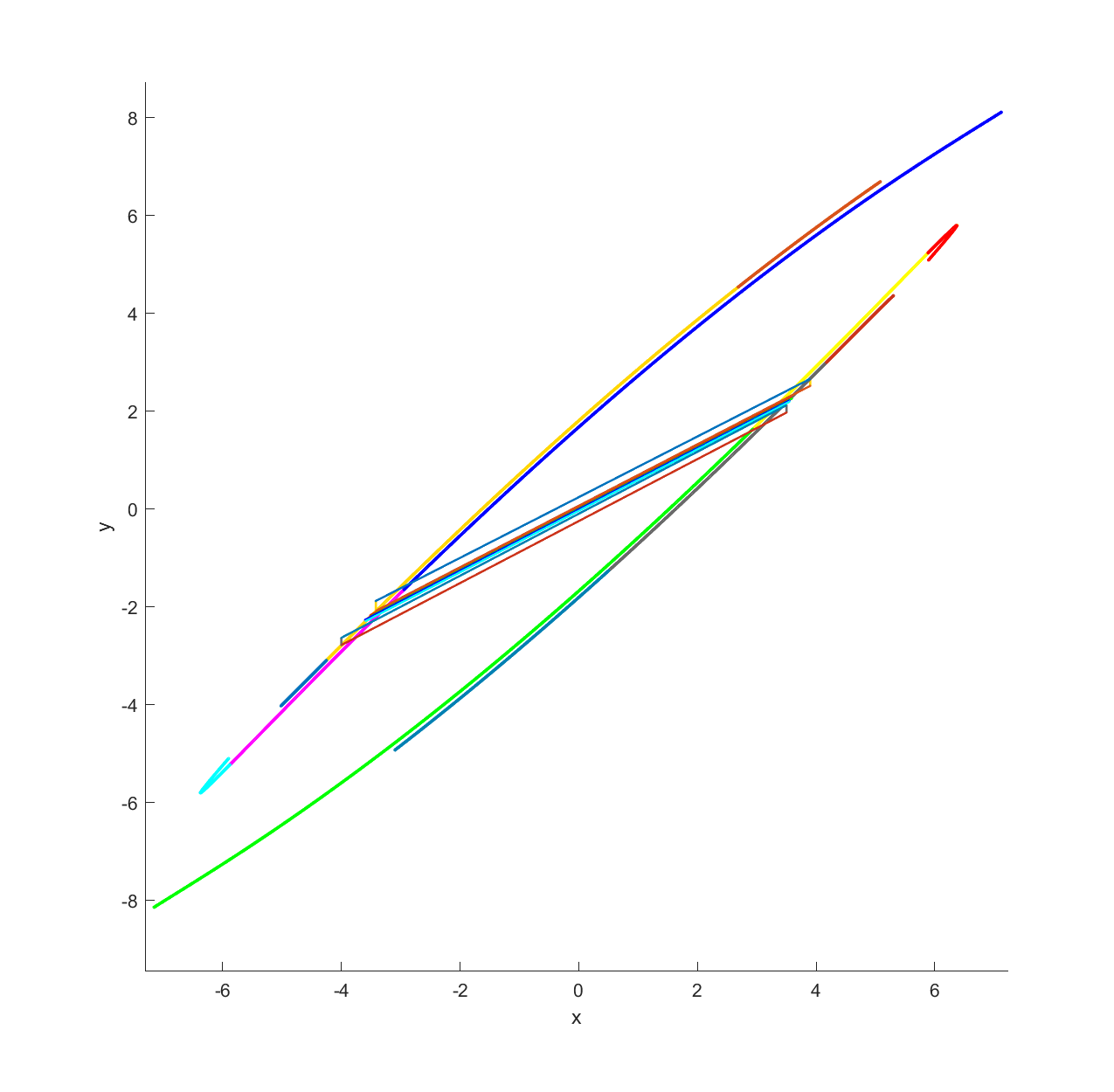}
	\caption{Four blocks and their images under the first return map.}
	\label{C4-1}
\end{figure}

\begin{figure}[h]
	\centering
	\includegraphics[height=8.5cm]{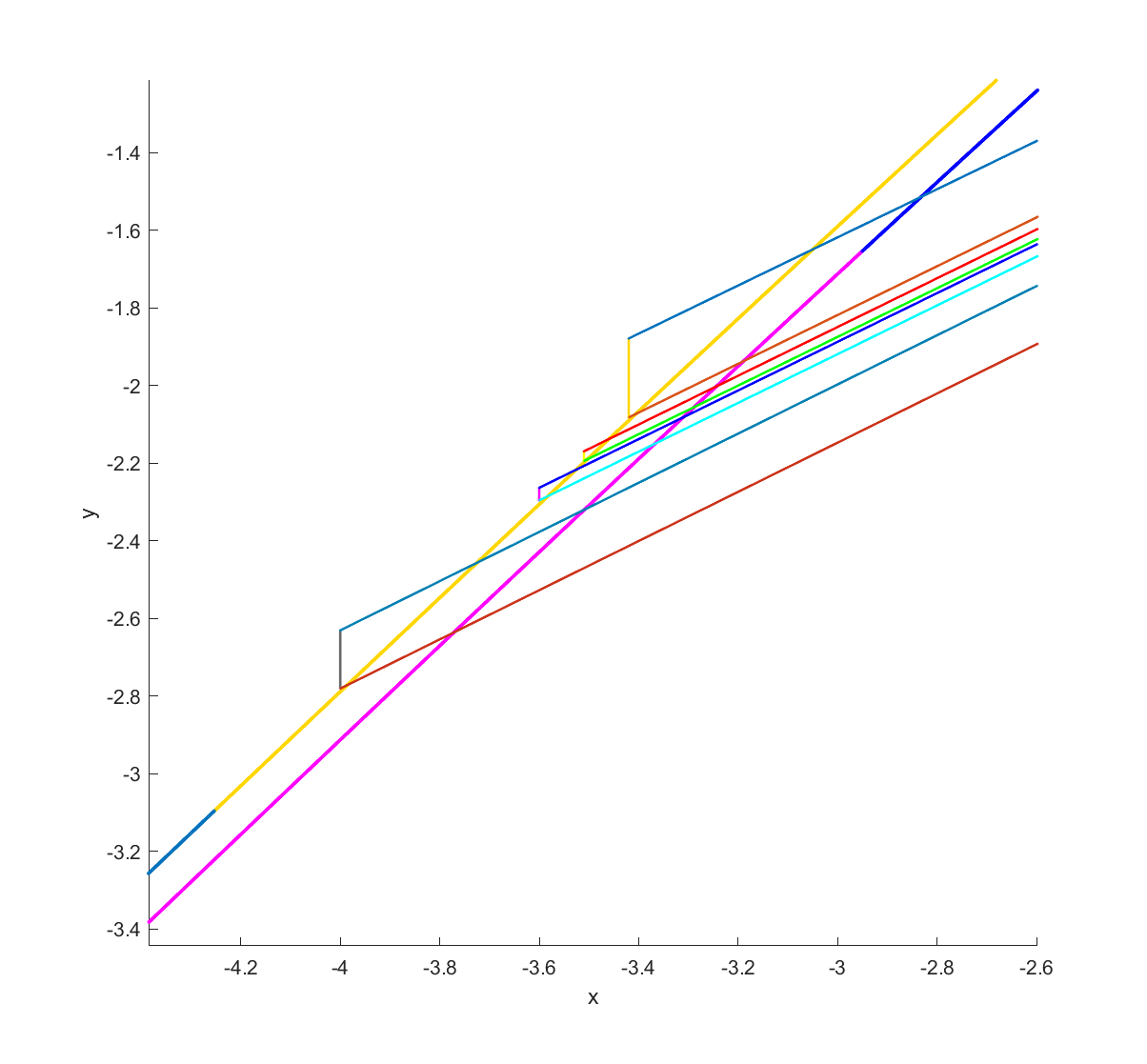}
	\caption{Magnified view of the left portion of Figure \ref{C4-1}.}
	\label{ppt-3.2}
\end{figure}

\begin{figure}[h]
	\centering
	\includegraphics[height=8.5cm]{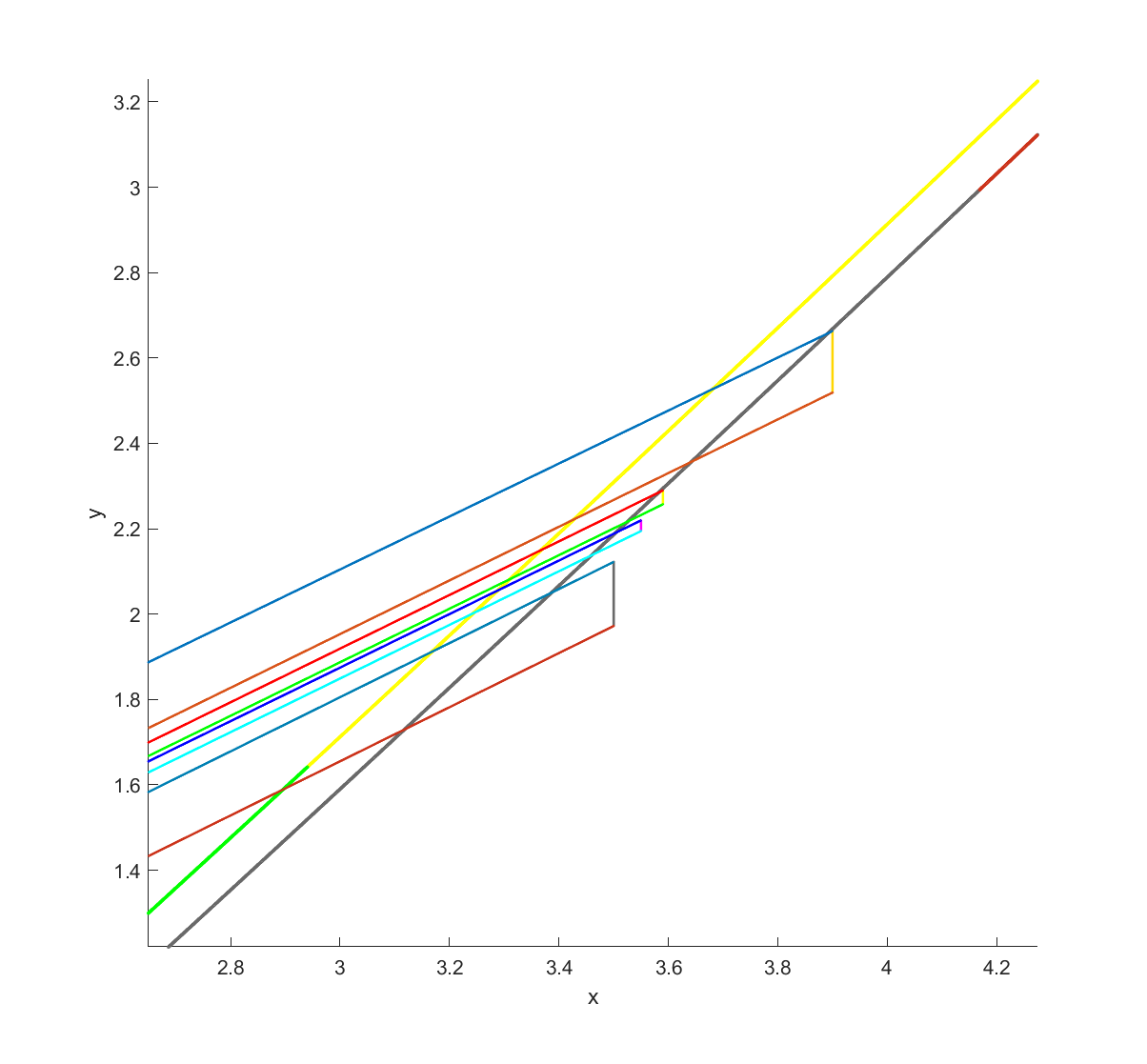}
	\caption{Magnified view of the right portion of Figure \ref{C4-1}.}
	\label{ppt-4.2}
\end{figure}

\begin{figure}[h]
	\centering
	\includegraphics[height=10cm]{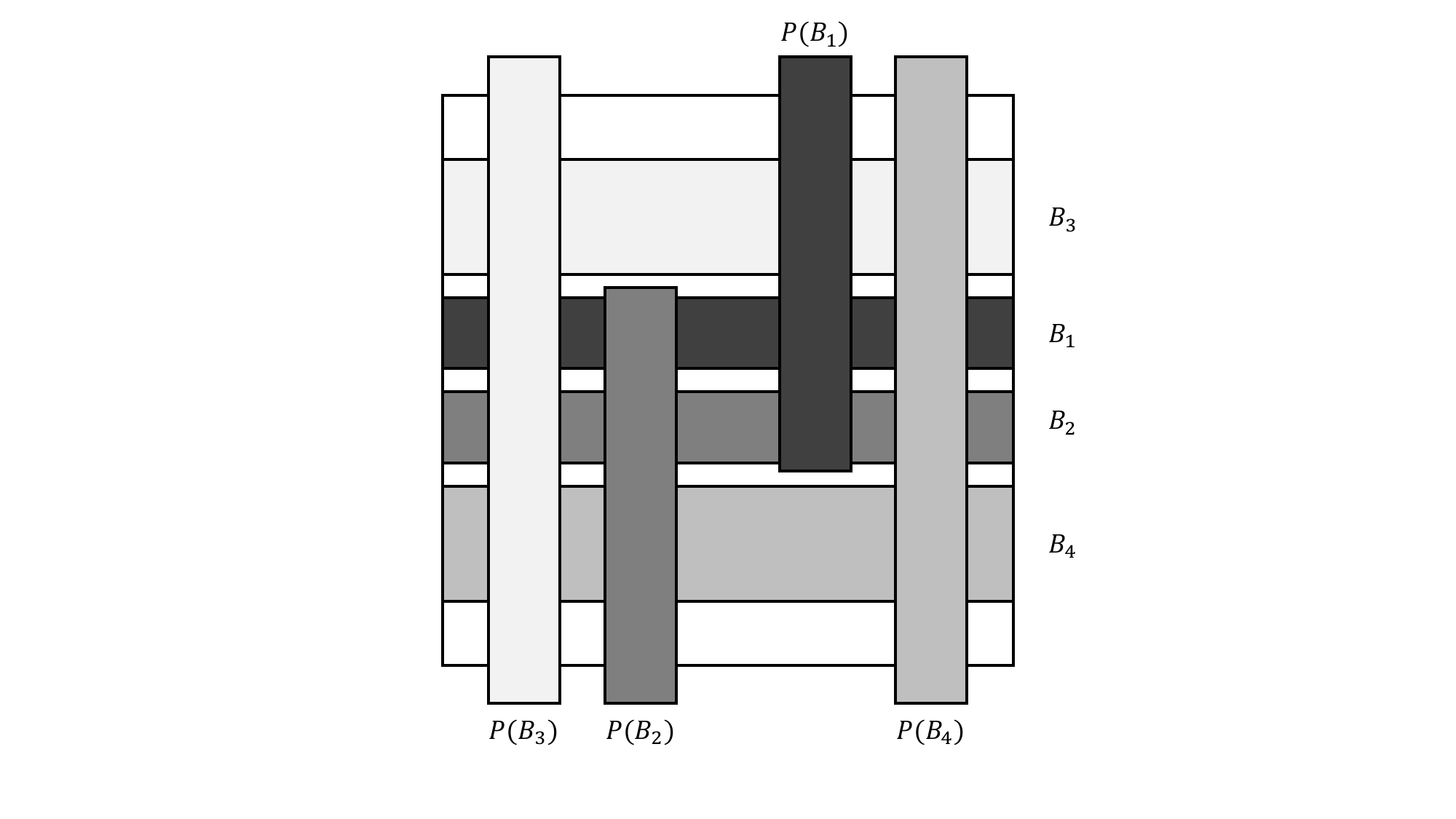}
	\caption{Illustration of incomplete crossing of four blocks in Figure \ref{C4-1}.}
	\label{ppt-5}
\end{figure}

As another example, we consider a polynomial system in \cite{Chen1999}, the so-called Chen system, represented by: 

\begin{eqnarray}
	\left\{\begin{array}{l}
		\dot{x}=35(y-x)\\
		\dot{y}=-7x+28y-xz\\
		\dot{z}=-3z+xy.
	\end{array}\right.
\end{eqnarray}
We select a horizontal Poincaré section $\{(x,y,21): -10\leq x\leq 10, -10\leq y\leq 10, xy<63\}$.
Four blocks and their images under the first return map are shown in Figures \ref{C4-1}--\ref{ppt-4.2}.
It is noteworthy that Figure \ref{ppt-3.2} and Figure \ref{ppt-4.2} imply that there is a semi-horseshoe involving four blocks in this polynomial system, as illustrated in Figure \ref{ppt-5}.
In \cite{Cheng2024}, the topological entropy of the Poincaré map \( P \) was estimated through iteration, yielding
\[h(P)\geq \frac12 \ln 6,\]
where $h(P)$ denotes the topological entropy of $P$.
With the theory presented in Section \ref{ITHT}, we can derive a more precise lower bound for \( h(P) \).


According to Theorem \ref{ITH}, $P$ is semi-conjugate to a subshift whose adjacency matrix is 
\begin{equation*}
	A=\left[
	\begin{array}{cccc}
1 & 1 & 1 & 0\\
1 & 1 & 0 & 1\\
1 & 1 & 1 & 1\\
1 & 1 & 1 & 1
	\end{array}
	\right]
\end{equation*}
The characteristic polynomial of $A$ is 
\[2 \lambda^2 - 4 \lambda^3 + \lambda^4,\]
which means the eigenvalues of $A$ are 
\[\lambda_1=\lambda_2=0, \lambda_3=2-\sqrt{2}, \lambda_4=2+\sqrt{2}.\]
Through further computation, we can obtain the spectral radius of $A$:
\[\rho(A)=\max\limits_{i}\|\lambda_i\|=2+\sqrt{2},\]
which implies that
\[h(P)\geq \ln (2+\sqrt{2}).\]

\begin{remark}
	Similar to the discussion in Remark \ref{remark1}, the result above can also be obtained by induction. When considering the iteration of $P$, we have the following table, where we set the number of blocks whose images under $P^n$ cross $B_1$ is $x_n$, and the number of blocks whose images under $P^n$ cross $B_3$ is $y_n$. Due to the symmetry of the structure, the number of blocks whose images under $P^n$ cross $B_2,B_4$ is $x_n,y_n$. Summary of iteration results are shown in Table \ref{table}.

\begin{table}
	\caption{\label{table}Summary of iteration results.}
	\begin{indented}
		\item[]\begin{tabular}{@{}lllll}
			\br
			Number of the blocks whose images crosses $B_i$ & $B_1$ & $B_2$ & $B_3$ & $B_4$   \\
			\mr
		$P$& 4 & 4 & 3 & 3 \\
		
		$P^2$& 14 & 14 & 10 & 10 \\
		
		\vdots& \vdots & \vdots & \vdots & \vdots \\
		
		$P^n$ & $x_n$ & $x_n$ & $y_n$ &$y_n$ \\
		
		$P^{n+1}$ & $2(x_n+y_n)$ & $2(x_n+y_n)$ & $x_n+2y_n$ &$x_n+2y_n$ \\
		\br 
		\end{tabular}
	\end{indented}
\end{table}

	So, we have
	\begin{eqnarray*}
		\left\{\begin{array}{l}
			x_{n+1}=2(x_n+y_n),\\
			y_{n+1}=x_n+2y_n,\\
			x_1=4,\\
			y_1=3.
		\end{array}\right.
	\end{eqnarray*}	
	By simple algebraic calculation, we have
	\begin{eqnarray*}
		\left\{\begin{array}{l}
			x_n=\frac{4-3\sqrt{2}}{2}(2-\sqrt{2})^{n-1}+\frac{4+3\sqrt{2}}{2}((2+\sqrt{2})^{n-1},\\
			y_n=\frac{3-2\sqrt{2}}{2}(2-\sqrt{2})^{n-1}+\frac{3+2\sqrt{2}}{2}(2+\sqrt{2})^{n-1}.
		\end{array}\right.
	\end{eqnarray*}	
	Similar as in Section \ref{TwoBlocks}, there will be $2y_n$ blocks whose images cross $B_i,(i=1,2,3,4)$ under $P^{n+1}$.
	Thus, the topological entropy of $P^n$ is at least $\ln(2y_{n-1})$, it follows that
	\begin{eqnarray*}
		h(P)&\geq \lim\limits_{n\to\infty}\frac1n\ln(2y_{n-1})\\
		&=\lim\limits_{n\to\infty}\frac1n\ln [(3-2\sqrt{2})(2-\sqrt{2})^{n-2}+(3+2\sqrt{2})(2+\sqrt{2})^{n-2}]\\
		&=\ln (2+\sqrt{2}).
	\end{eqnarray*}
\end{remark}

\section{Summary}\label{summary}

In this paper, we have generalized the topological horseshoe theory. We introduce the concepts of incomplete crossing and semi-topological horseshoe and develop a theory of semi-topological horseshoes by virtue of the theory of subshifts of finite type in symbolic dynamical systems.
As applications of our theory, we revisit the perturbed Duffing system and a polynomial system, providing more accurate estimates of the lower bounds of the topological entropy of their Poincaré maps. The presence of incomplete crossing structures in these systems suggests that such semi-horseshoes are prevalent in chaotic systems. We hope this work will stimulate further investigations about chaos.


\section*{References}

\begin{thebibliography}{16}  
	

	\bibitem[1]{Adler1965}	Adler R L, Konheim A G and McAndrew M H 1965 Topological entropy {\it Trans. Amer. Math. Soc.} {\bf 114} 309-319
	
	
	
	
	
	\bibitem[2]{Chen1999} Chen G R  and Ueta T 1999 Yet another chaotic attractor {\it Internat. J. Bifur. Chaos Appl. Sci. Engrg.} {\bf 9} 1465--1466
	
	\bibitem[3]{Cheng2024} Cheng J and Yang  X-S 2024 New Topological Horseshoe in the Chen System {\it Internat. J. Bifur. Chaos Appl. Sci. Engrg.}  2430018
	
	\bibitem[4]{Douglas2021} Douglas L and Brian M 2021 {\it An Introduction to Symbolic Dynamics and Coding} (Cambridge: Cambridge University Press)
	
	
	
	\bibitem[5]{Guckenheimer1984} Guckenheimer J and Holmes P 1984 {\it Nonlinear oscillations, dynamical systems, and bifurcation of vector fields} (Applied Mathematical Sciences)
	
	
	
	\bibitem[6]{Kennedy2001-1} Kennedy J and Yorke J A 2001 Topological horseshoes {\it Trans. Amer. Math. Soc.} {\bf 353}, pp.  2513--2530
	
	\bibitem[7]{Kennedy2001-2} Kennedy J, Koçak S and Yorke J A 2001 A chaos lemma {\it Amer. Math. Monthly} {\bf 108} 411--422
	
	\bibitem[8]{Lorenz1963} Lorenz E N  1963  Deterministic nonperiodic flow {\it J. Atmos. Sci.} {\bf20}   130--141
	
	\bibitem[9]{Melnikov1963} Melnikov  V K 1963 On the stability of a center for time-periodic perturbations {\it Tr. Mosk. Mat. Obs.}{\bf 12},   3--52
	
	
	\bibitem[10]{Smale1965} Smale  S 1965  Diffeomorphisms with Many Periodic Points {\it Differential and Combinatorial Topology (A Symposium in Honor of Marston Morse)} (Princeton: Princeton University Press)  63--80
	
	
	\bibitem[11]{Smale1967} Smale  S 1967  Differentiable dynamical systems {\it Bull. Amer. Math. Soc.} {\bf 73} 747-817
	
	\bibitem[12]{Tucker2002}Tucker  W 2002 A rigorous ODE solver and Smale's 14th problem {\it Found. Comput. Math.} {\bf2} 53--117
	
	\bibitem[13]{Wiggins1988} Wiggins S 1988 {\it Global Bifurcations and Chaos:Analytical Methods} {\it Appl. Math. Sciences} vol 73 ( New York : Springer-Verlag)
	
	\bibitem[14]{Yang2004} Yang X-S and Tang Y [2004]  Horseshoes in piecewise continuous maps  {\it Chaos Solitons Fractals} {\bf 19} 841--845
	
	\bibitem[15]{Yang2009} Yang X-S 2009 Topological horseshoes and computer assisted verification of chaotic dynamics  {\it Internat. J. Bifur. Chaos Appl. Sci. Engrg.} {\bf 19} 1127--1145
	
	\bibitem[16]{Zhou1997} Zhou  Z   1997  {\it Symbolic Dynamics} (Shanghai: Shanghai Scientific and Technological Education Publishing House)
	
	
\end{thebibliography}
\end{document}